\documentclass[12pt]{report}
\usepackage{amsfonts,amssymb,amsmath}
\usepackage[T2A]{fontenc}  
\usepackage[cp1251]{inputenc}
\usepackage[english,russian]{babel}
\usepackage{geometry}

\newcommand{\RR}{{\rm I\hskip -.15em R}}
\newfont{\ruk}{eusm10 at 12pt}

\def\Div{{\rm div}}

\def\O{\Omega}
\def\beq{\begin{equation}}
\def\eeq{\end{equation}}

\def\nab{\nabla}
\def\f{\frac}

\def\ve{\varepsilon}
\def\vp{\varphi}

\def\ilo{\int\limits_{\O}}

\def\ds{\displaystyle}

\def\l{\left}
\def\r{\right}

\def\C0{C_0^\infty}

\def\pa{\partial}

\newtheorem{theorem}{Teoрeма}[section]
\newtheorem{lemma}{Лемма}[section]

\newcommand{\doc}{{\bf Доказательство. }}

\begin{document}

\begin{center}
\large\textbf{On homogenization estimates in Neuman boundary value problem\\ for an elliptic equation with multiscale coefficients}
\end{center}
\medskip
\begin{center}
{\sc{    S.~E.~Pastukhova, R.~N.~Tikhomirov}}
\end{center}
\medskip

Homogenization of a scalar elliptic equation in a bounded domain with Neuman boundary condition is studied. Coefficients of the operator are oscillating over two  different groups of variables with different small periods $\ve$ and $\delta=\delta(\ve)$. We assume that ${\delta}/{\ve}$ tends to zero as $\ve$ tends to zero. It is known that the limit problem is obtained through reiterated homogenization procedure and corresponds to an elliptic equation with constant coefficients. The difference for resolvents of the initial  and the  limit problems is estimated in  operator $(L^2\to L^2)$-norm. This estimate is of order $\tau=\tau(\ve)$ which is the maximum of $\ve$ and ${\delta}/{\ve}$. We find also the approximation of the initial resolvent in operator $(L^2\to H^1)$-norm, it is of order $\sqrt{\tau}$.

%
%
%
%

\bigskip  

В теории усреднения наблюдается повышенный интерес к $L^2$- и $H^1$-оценкам  операторного типа для задач с осциллирующими коэффициентами, начавшийся после выхода работы М.Ш. Бирмана и Т.А. Суслиной \cite{BirSu}. Оценки называются операторными, потому что их можно сформулировать  в операторных нормах в терминах резольвент исходной и усредненной задач и, быть может, некоторых корректирующих операторов.
Настоящая работа продолжает цикл работ \cite{Zh1}-\cite{PR}, где 
операторные оценки усреднения доказаны для разнообразных уравнений модифицированным методом первого приближения, предложенным В.В. Жиковым в \cite{Zh1}. Отличительные черты этого метода -- введение дополнительного параметра интегрирования (за счет прямого сдвига или в скрытом виде  за счет сглаживания по Стеклову), а также специальный анализ первого приближения, или  $H^1$-приближения.
В указанных работах 
были изучены эллиптические уравнения резольвентного типа во всем пространстве, в том числе, векторные уравнения теории упругости, а также соответствующие краевые задачи в области с условием Дирихле и Неймана на границе. В данной статье рассматривается  скалярное эллиптическое уравнение с двухмасштабными коэффициентами. Коэффициенты уравнения  осциллируют по двум группам переменных с периодами $\ve$ и $\delta=\delta(\ve)$ разного порядка малости при $\ve\to 0$. Будем считать для определенности $\delta(\ve)$ большего порядка малости, чем $\ve$, то есть $\lim\limits_{\ve\to 0}\f{\delta(\ve)}{\ve}=0$. Для задачи Неймана в ограниченной области  получена $L^2$-оценка погрешности усреднения. Эта оценка порядка $\max\{\ve,\f{\delta}{\ve}\}$. Она является следствием  $H^1$-оценки  погрешности усреднения, которая имеет порядок $\max\{\sqrt{\ve},\sqrt{\f{\delta}{\ve}}\}$.
Чтобы доказать  оценки,
используется $H^1$-приближение, в котором участвует  сглаживание по Стеклову.
Это приближение
 для скалярной задачи  может быть выбрано  
 более простого вида, чем в общем случае, охватывающем и векторные задачи. Хотя задача двухмасштабна, а именно, содержит два типа быстрых переменных, сглаживание по Стеклову берется только 
 с одним параметром сглаживания. 
  Такое сглаживание  
 сильно упрощает конструкцию $H^1$-приближения.

\section{Постановка задачи и основной результат} Введем пространства функций с условием ортогональности
$$
\tilde{L}^2(\O)=\{u\in {L}^2(\O):\, \ilo u\,dx=0\},\quad
\tilde{H}^1(\O)=\{u\in {H}^1(\O):\, \ilo u\,dx=0\}.
$$
Рассмотрим  краевую задачу Неймана в ограниченной липшицевой области $\Omega\in \RR^d$  для эллиптического уравнения второго порядка
\beq\label{eq1}
u^\ve\in \tilde{H}^1(\O),\quad
\ilo a^\ve \nab u^\ve\nab  \varphi\,dx=\ilo f\varphi\,dx\quad \forall \varphi\in H^1(\O),\quad
f\in \tilde{L}^2(\O).
\eeq
Матрица $a^\ve(x)$ не обязательно симметрична и  имеет "двухмасштабную"  структуру
\beq\label{eq2}
a^\ve(x)=a\left(\frac{x}{\ve},\frac{x}{\delta}\right),
\eeq
где элементы матрицы $a(y,z)$ измеримы и периодичны по $y$ и $z$,  ячейкой периодичности служит куб $Y=Z=[-\frac12,\frac12)^d$. При $\ve\to 0$ и $\delta\to 0$ в задаче (\ref{eq1}) имеем быстро осциллирующие коэффициенты, $\ve$- и $\delta$-периодические по первой и второй группе переменных. Для определенности
$$
\delta=\delta (\ve),\quad \lim\limits_{\ve\to 0}\frac{\delta(\ve)}{\ve}=0.
$$
Дополним условия на $a(y,z)$:
\beq\label{eq3}
\ds{|a(y,z)-a(y',z)|\le c_L |y-y'|,}\atop
\ds{\mu|\xi|^2\le a(y,z)\xi\cdot\xi,\quad a(y,z)\xi\cdot\eta\le \mu^{-1}|\xi||\eta|\ \forall\xi,\eta\in\RR^d,\ \mu>0,}
\eeq
для почти всех $z\in \RR^d$ и всех $y,y'\in \RR^d$.

Существование единственного решения задачи (\ref{eq1}) гарантировано условиями ортогональности  и получается по теореме Лакса -- Мильграма. Выполнено энергетическое равенство
$$
\ilo  a^\ve \nabla u^\ve\nabla  u^\ve\,dx=\ilo fu^\ve\,dx,
$$
из которого следует энергетическая оценка
$$
\|u^\ve\|_{H^1(\O)}\le c\|f\|_{L^2(\O)},
$$
в силу  неравенства Пуанкаре
$$
\ilo  |u|^2\,dx\le C_P\l(\ilo  |\nabla u|^2\,dx+\l(\ilo  u\,dx\r)^2\r).
$$

С задачей (\ref{eq1}) связана  усредненная задача
\beq\label{eq4}
u\in \tilde{H}^1(\O),\quad
\ilo a^0 \nabla u^0\nabla  \varphi\,dx=\ilo f\varphi\,dx\quad \forall \varphi\in H^1(\O),\quad
f\in \tilde{L}^2(\O),
\eeq
где  $a^0$ -- постоянная эллиптическая  матрица, ее точное определение дано ниже через решения вспомогательных задач на ячейке периодичности.

Будем предполагать, что область достаточно гладкая, например класса $C^{1,1}$, так что по эллиптической теории решение задачи (\ref{eq4}) принадлежит пространству $H^2(\O)$ с оценкой
\beq\label{eq5}
\|u\|_{H^2(\O)}\le c_0\|f\|_{L^2(\O)},\ c_0=const (d,\mu,\O),
\eeq
кроме того, существует линейный ограниченный оператор продолжения
\beq\label{eq6}
\mathcal{P}:\ H^2(\O)\to H^2(\RR^d).
\eeq
Всегда считаем, что решение задачи (\ref{eq4}), если нужно, продолжено во внешность области $\O$ с помощью оператора $\mathcal{P}$ и это продолжение также обозначаем через $u$.

Справедлива
\begin{theorem}
Для разности решений (\ref{eq1}) и (\ref{eq4}) имеет место оценка
\beq\label{eq7}
\|u^\ve-u\|_{L^2(\O)}\le C\max\{\ve,\frac{\delta}{\ve}\}\|f\|_{L^2(\O)},
\eeq
где константа $C$ зависит лишь от размерности $d$, постоянных $\mu$,  $c_L$ и области $\O$.
\end{theorem}
Оценка (\ref{eq7}) выводится из соответствующей $H^1$-оценки (см. (\ref{eq23})), доказательству которой посвящена  основная часть настоящей работы.

\section{Задачи на ячейке}
 Матрица $a^0$ определяется в процессе повторного усреднения. Введем некоторые обозначения: $e^1,\ldots,e^d$ -- канонический базис в $\RR^d$; $\langle\cdot\rangle_Z=\int\limits_Z\cdot dz$ -- среднее по ячейке; $H^1(Y)=H^1_{per}(Y)$, $H^1(Z)= H^1_{per}(Z)$, $H^2(Y)=H^2_{per}(Y)$ и т.д. -- соболевские пространства периодических функций.

Нам понадобятся следующие задачи на ячейке периодичности:
\beq\label{eq8}
M_j(y,\cdot)\in H^1(Z),\ \Div_z [a(y,z)(e^j+\nabla_z M_j(y,z))]=0,\ \langle M_j (y,\cdot)\rangle_Z=0,
\eeq
\beq\label{eq9}
N_j\in H^1(Y),\ \Div_y [\hat{a}(y)(e^j+\nabla N_j(y))]=0,\ \langle N_j \rangle_Y=0, \atop j=1,2,\ldots,d.
\eeq
Здесь
$$
\hat{a}(y)=\langle a(y,\cdot)(I+\nabla_z M(y,\cdot))\rangle_Z
$$
является "промежуточной"\ усредненной матрицей,  $I$ -- единичная $(d\times d)$-матрица, а $M(y,z)$ -- вектор с координатами $M_j(y,z)$. Окончательной усредненной матрицей будет
$$
a^0=\langle \hat{a}(I+\nabla_y N)\rangle_Y.
$$
Усредненные матрицы $\hat{a}$ и $a^0$ подчинены условиям типа (\ref{eq3}), так что обе задачи на ячейке, а также усредненная задача поставлены корректно.

Решения задач  на ячейке 
 обладают следующими свойствами.

\begin{lemma}\label{l1}
Пусть $M_j(y,z)$ и $N_j(y)$ -- решения уравнений (\ref{eq8}) и (\ref{eq9}). Тогда\\
$i)\ \|M_j(y,\cdot)\|_{H^1(Z)}\le c$, $c=const(\mu)$;\\
$ii)\ M_j(y,\cdot)$ -- липшицева по $y\in Y$ функция со значениями в $H^1(Z)$ с константой Липшица, зависящей лишь от постоянных $\mu$ и $c_L$;\\
$iii)$ матрица $\hat{a}(y)$ -- липшицева  с константой $c=const(\mu,c_L)$;\\
$iv)\ N_j,\ \nabla_y N_j\in L^\infty (Y)$;\\
$v)\  \nabla^2_y N_j\in L^p(Y)\ \forall p\ge 1$;\\
$vi)\ \| M_j(y,\cdot)\|_{L^\infty (Z)}\le c$ для всех $y$.
\end{lemma}
Доказательство леммы приведено в \S\!! 7.

Введем периодические векторы, связанные с задачами на ячейке,
\beq\label{eq10}
p^j (y,z)=a(y,z)(e^j+\nabla_z M_j(y,z))-\hat{a}(y) e^j,\ \quad g^j(y)=\hat{a}(y)(e^j+\nabla_y N_j (y)) - a^0e^j.
\eeq
Первый из них представляет собой  соленоидальный по $z$ вектор с нулевым средним:
\beq\label{eq11}
\langle p^j(y,\cdot)\rangle_Z =0,\quad  \Div_z p^j(y,z)=0
\eeq
для всех $y$,
второй -- также  соленоидальный  вектор с нулевым средним:
\beq\label{eq12}
\langle g^j\rangle_Y =0,\ \Div g^j(y)=0.
\eeq
Известно \cite[гл. I, \S 1]{ZhKO}, что такие векторы можно представить в виде дивергенции от кососимметрической матрицы.
В наших предположениях это представление, для векторов из (\ref{eq10}), обладает дополнительными свойствами.

\begin{lemma}
Пусть $g^j$  и $p^j$ -- векторы из (\ref{eq10}).  Тогда найдутся такие кососимметрические матрицы $G^j\in H^1 (Y)^{d\times d}$ и $P^j(y,\cdot)\in H^1(Z)^{d\times d}$, что справедливы представления:\\
$
a)\quad g^j(y)=\Div\, G^j(y);\\
$
$b)\quad p^j(y,z)=\Div_z P^j(y,z)$.\\
При этом
\beq\label{eq13}
\|G^j\|_{H^1(Y)^{d\times d}}\le c\|g^j\|_{L^2(Y)^d},\atop
\|P^j(y,\cdot)\|_{H^1(Z)^{d\times d}}\le c \|p^j(y,\cdot)\|_{L^2(Z)^d}
\eeq
с константой $c=const (d,\mu)$. Кроме того,
\beq\label{eq1.3.3.18}
g^j\in L^\infty (Y)^d,\ G^j\in L^\infty (Y)^{d\times d}
\eeq
и матрица $P^j(y,z)$  липшицева по $y$ со значениями в $H^1(Z)^{d\times d}$, что влечет  оценку
для ее элементов
\beq\label{eq14}
\|\nabla_y P^j_{ik}(y,\cdot)\|_{H^1(Z)^d}\le c,\quad c=c(\mu,c_L),
\eeq
для почти всех $y$.
\end{lemma}
\doc
Сопоставим вектору $g^j=\{g^j_{k}\}$
 задачу на ячейке $Y$
\beq\label{eq150}
\Delta\varphi^j_k=g^j_k,\ \langle\varphi^j_k\rangle_Y=0,\ k=1,\ldots,d.
\eeq
Решение $\vp^j_k$ существует и единственно. Справедлива оценка
\beq\label{eq15}
\|\varphi^j_k\|_{H^2(Y)}\le c\|g^j_k\|_{L^2(Y)}.
\eeq
В силу соленоидальности вектора $g^j$ выполнено соотношение
$$
\Delta\frac{\partial\varphi^j_k}{\partial y_k}=\frac{\partial g^j_k}{\partial y_k}=0, \ k=1,\ldots d,
$$
(как обычно, по повторяющимся индексам подразумевается суммирование от $1$ до $d$, если не оговорено противное).
 Таким образом, $\Div\vp^j$ -- периодическая гармоническая функция,  $\langle\Div\vp^j\rangle_Y=0$, следовательно,  $\Div\vp^j=0$.

Введем 
кососимметрическую матрицу $G^j=\{G^j_{ik}\}$, такую что
\beq\label{eq16}
G^j_{ik}=\frac{\partial\varphi^j_k}{\partial y_i}-\frac{\partial\varphi^j_i}{\partial y_k},
\eeq
и проверим равенство $\Div G^j=g^j$. В самом деле,
$$
\frac{\partial G^j_{ik}}{\partial y_i}=\Delta\varphi^j_k-\frac{\partial}{\partial y_k}\Div\varphi^j=\Delta\varphi^j_k=g^j_k,
$$
и представление a) получено. Оценка (\ref{eq13})$_1$ следует из (\ref{eq16}), (\ref{eq15}).

Свойство (\ref{eq1.3.3.18})$_1$ следует из опеределения $g^j$
и ограниченности $\nabla N_j$ (см. лемму 2.1,  $iv)$ ). Тогда свойство (\ref{eq1.3.3.18})$_2$ 
обеспечено  эллиптической теорией и теоремами вложения: как решение уравнения (\ref{eq150})
функция $\varphi^j_k$ принадлежит $W^{2,p}(Y),\, \forall p\ge 1,$ и, значит,  ограничены
ее производные, а также  $G^j_{ik}$  в силу  (\ref{eq16}).

Аналогично, матрица $P^j =\{P^j_{ik}\} $ определяется своими элементами
$$
P^j_{ik}=\frac{\partial\psi^j_k}{\partial z_i}-\frac{\partial\psi^j_i}{\partial z_k},
$$
где $\psi^j_{k}(y,\cdot) \in H^2(Z)$ есть решение уравнения
\beq\label{eq17}
\Delta_z\psi^j_k (y,z)=p^j_k(y,z),\ \langle\psi^j_k(y,\cdot)\rangle_Z=0,\ k=1,\ldots,d.
\eeq
Правая часть $\psi^j_k(y,\cdot)$ липшицева по $y$ со значениями в $L^2(Z)$, согласно определению (\ref{eq10})$_1$ и свойствам липшицевости участвующих в нем функций $a(y,z)$, $M(y,z)$ и $\hat{a}(y)$. Из (\ref{eq17}) вытекает уравнение
$$
\Delta_z[\psi^j_k(y+h,z)-\psi^j_k(y,z)] = p^j_k(y+h,z)-p^j_k(y,z),
$$
из которого получаем оценки
$$
\|\psi^j_k(y+h,\cdot)-\psi^j_k(y,\cdot)\|_{H^2(Z)}\le c\|p^j(y+h,\cdot)-p^j(y,\cdot)\|_{L^2(Z)^d}\le C|h|,
$$
$$
\|P^j(y+h,\cdot)-P^j(y,\cdot)\|_{H^1(Z)^{d\times d}}\le c\|p^j(y+h,\cdot)-p^j(y,\cdot)\|_{L^2(Z)^d}\le C|h|,
$$
что доказывает липшицевость $P^j(y,\cdot)$ со значениями в $H^1 (Z)^{d\times d}$.  
Следовательно, градиент $\nabla_y P^j_{ik}(y,\cdot)$  существует как ограниченная функция со значениями в $H^1 (Z)^{d}$ и справедливо соотношение (\ref{eq14}).
Лемма доказана.

\section{Первое приближение}

В соответствии с оценкой (\ref{eq7}) функцию $u$ называют 
$L^2$-приближением к $u^\ve$. Чтобы аппроксимировать $u^\ve$ в $H^1$-норме, требуется к нулевому приближению $u$ добавить корректор. Так получаем первое приближение
\beq\label{eq18}
v^\ve(x)= u(x)+\ve N(y)\cdot\nabla u(x)+\delta M(y,z)\cdot (I+\nabla_y N(y))\nabla u(x),\quad y=\f{x}{\ve},\ z=\f{x}{\delta},
\eeq
которое рассматривалось, например в \cite{BLP}. Посмотрим, насколько регулярна функция $v^\ve$ в наших предположениях.

Свойства    функций $N_j$,  $M_j$ (см. лемму 2.1, iv), vi)) обеспечивают принадлежность $v^\ve$ пространству $L^2(\O)$. Действительно, для составляющих корректора имеют место оценки
\beq\label{eq19}
\ilo \l|N(\f{x}{\ve})\cdot\nabla u(x)\r|^2dx\le c\ilo|\nabla u(x)|^2dx,
\eeq
$$
\ilo \l|M(\f{x}{\ve},\f{x}{\delta})\cdot (I+\nabla_y N(\f{x}{\ve}))\nabla u(x)\r|^2dx\le c \ilo |\nabla u(x)|^2 dx,
$$
так как $N$, $\nabla_y N_j\in L^\infty (Y)^d$ и   $M(y,\cdot)\in L^\infty (Z)$ для всех $y$. Для градиента $\nabla v^\ve$ подобная оценка в $L^2(\O)$ не пройдет из-за вхождения в его структуру неограниченных функций  $\nabla_y M$, $\nabla_z M$ и $\nabla^2 N$. Чтобы обойти это препятствие будем брать вместо (\ref{eq18}) первое приближение со сглаженным корректором, а именно:
\beq\label{eq20}\ds{
\hat{v}^\ve (x)=u(x)+ \ve \int\limits_{Z} N(y')\cdot \nabla u(x') d\sigma +\delta \int\limits_{Z} M(y',z)\cdot (I+\nabla_y N(y'))\nabla u(x') d\sigma,}
\atop\ds{
x'=x-\delta\sigma,\ y'=\frac{x}{\ve}-\f{\delta}{\ve}\sigma,\ z=\frac{x}{\delta}.}
\eeq
Будем также использовать более компактное представление $\hat{v}^\ve(x)$ в виде суммы  нулевого приближения $u(x)$ и корректора $K_\ve (x)$:
\beq\label{eq1.3.3.44}
\hat{v}^\ve (x) = u(x) + K_\ve (x),\ K_\ve (x) = \ve K_{\ve,1}(x) +\delta(\ve) K_{\ve,2}(x).
\eeq

Объясним роль дополнительного параметра интегрирования, введенного в первое приближение  (\ref{eq20}).
\begin{lemma}
Пусть  $b(y,z)\in L^\infty (Y, L^2_{per}(Z))$ и $u(x)\in L^2(\RR^d)$. Тогда функция
$$
w(x)=\int\limits_Z b(y',z) u(x') d\sigma,\quad y'=\f{x}{\ve}-\f{\delta}{\ve}\sigma,\ z=\f{x}{\ve},\ x'=x-\delta\sigma,
$$
принадлежит $L^2(Q)$ с оценкой
\beq\label{eq21}
\|w\|_{L^2(Q)} \le \sup\limits_y \|b(y,\cdot)\|_{L^2(Z)}\|u\|_{L^2(Q_{\ve})}.
\eeq
где $Q$ -- произвольная область в $\RR^d$, а $Q_{\ve}$ --  $\ve$-окрестность $Q$.
\end{lemma}
\doc 
Используя неравенство Коши -- Буняковского  и замену переменных $x\to x'=x-\delta\sigma$, получим цепочку соотношений 
\beq\label{eq210}
\ds{\int\limits_{Q} |w(x)|^2 dx\le\int\limits_{Q}\int\limits_{Z} \l|b\l(\f{x}{\ve}-\f{\delta}{\ve}\sigma,\f{x}{\delta}\r)\r|^2 |u(x-\delta\sigma)|^2 dxd\sigma\le}
\atop\ds{
\le \int\limits_{Q_{\ve}}\int\limits_{Z} \l|b\l(\f{x}{\ve},\f{x}{\delta}+\sigma\r)\r|^2 |u(x)|^2 dxd\sigma =
\int\limits_{Q_{\ve}} |u(x)|^2\l(\int\limits_Z\l|b\l(\f{x}{\ve},\f{x}{\delta}+\sigma\r)\r|^2 d\sigma\r) dx.}
\eeq
Ввиду периодичности функции $b(y,\cdot)$, внутренний интеграл оценивается величиной $\sup\limits_y \|b(y,\cdot)\|^2_{L^2(Z)}$, которая выносится за пределы внешнего интеграла,  и в результате приходим к неравенству (\ref{eq21}). Лемма доказана.

Теперь можно показать, что $\nabla\hat{v}^\ve\in L^2(\O)$. Запишем градиент как сумму
$$
\nabla{v}^\ve(x)= \nabla u(x) +\int\limits_Z \nabla_y N(y')\nabla u(x')d\sigma+\int\limits_Z \nabla_z M(y',z)(I+\nabla_y N(y'))\nabla u(x')d\sigma+
$$
$$
+\f{\delta}{\ve}\int\limits_Z\nabla_y\l[M(y',z)\cdot(I+\nabla_y N(y'))\nabla u(x')\r]d\sigma+ \ve\int\limits_{Z} \nabla^2 u(x')N(y')d\sigma+
$$
$$
+\delta\int\limits_Z (I+\nabla_y N(y'))\nabla^2 u(x')M(y',z)d\sigma,\ \quad y'=\f{x}{\ve}-\f{\delta}{\ve}\sigma,\ x'=x-\delta\sigma,\ z=\f{x}{\delta},
$$
и изучим каждое слагаемое. Для второго,  пятого и  шестого слагаемых достаточно воспользоваться свойствами iv), vi)  из леммы 2.1 и оценкой (\ref{eq5}), а для третьего -- применить лемму 3.1. Четвертое слагаемое требует более детального анализа. Оно разбивается в сумму:
\beq\label{eq22}
\ds{
\int\limits_Z\nabla_y\l[M(y',z)\cdot(I+\nabla_y N(y'))\nabla u(x')\r]d\sigma =}
\atop\ds{
=\int\limits_Z\nabla_y M(y',z)(I+\nabla_y N(y'))\nabla u(x')d\sigma + \int\limits_Z \nabla^2_y N(y')\nabla u(x')\cdot M(y',z)d\sigma = T_1+T_2.}
\eeq
Используя неравенство Коши -- Буняковского  и  
ограниченность $\nabla N$,  имеем оценку
$$
|T_1|^2  \le c\int\limits_Z |\nabla_y M(y',z)\nabla u(x')|^2 d\sigma,
$$
к мажоранте которой применимо неравенство (\ref{eq210}). Слагаемое $T_2$
 оценивается несколько иначе. Прежде всего запишем неравенства
$$
\ilo|T_2|^2\,dx=
\ilo
|\int\limits_Z \nabla^2_y N(y')\nabla u(x')\cdot M(y',z)d\sigma
|^2 dx\le C\ilo\int\limits_Z |\nabla^2_y N(y')\nabla u(x')|^2d\sigma dx\le
$$
$$
\le C \int\limits_{\Omega_{\ve}}|\nabla^2_y N(y)|^2\l|\nabla u(x)\r|^2  dx\le C \int\limits_{\Omega_1}|\nabla^2_y N(y)|^2\l|\nabla u(x)\r|^2 dx,
$$
где использовали неравенство Коши -- Буняковского и замену переменных, а также учли ограниченность функции $N$ и перешли к области интегрирования $\Omega_1\supset\Omega_{ \ve}$.
Проверим, что последний интеграл ограничен. Перепишем его  в виде
$$
\int\limits_{\Omega_1} \beta(\f{x}{\ve})\Phi(x) dx.
$$
Здесь $\Phi\in  L^{1+\gamma}(\O_1)$, $\gamma >0$, по теореме вложения, так как $\nabla u\in H^1(\O_1)$,
 а  $\beta\in L^q_{per}(Y)$ при любом сколь угодно большом $q$ (см.  лемма 2.1, v)). Тогда по неравенству Гельдера
\beq\label{eq230}
\int\limits_{\Omega_1} \beta(\f{x}{\ve})\Phi(x) dx\le \l(\int\limits_{\Omega_1} |\beta(\f{x}{\ve})|^{q} dx\r)^{\f 1q}\l(\int\limits_{\Omega_1}|\Phi(x)|^{1+\gamma}dx\r)^{\f{1}{1+\gamma}} = C_{\Phi}\l(\int\limits_{\Omega_1} |\beta(\f{x}{\ve})|^{q} dx\r)^{\f 1q} ,\quad q=\f{1+\gamma}{\gamma},
\eeq
где
$$
\lim\limits_{\ve\to 0}\int\limits_{\Omega_1} \l|\beta\l(\f{x}{\ve}\r)\r|^{q} dx = |\O_1|\int\limits_{Y} \l|\beta(y)\r|^q dy
$$
по свойству среднего значения периодической функции. Величина  $C_{\Phi}=\|\Phi\|_{L^{1+\gamma}(\O_1)}$ контролируется нормой $\|\nabla^2 u\|_{L^2(\O_1)}$, где $u$ -- продолженное на $\RR^d$ с помощью оператора $\mathcal{P}$ (см. (\ref{eq6})) решение усредненной задачи.
В итоге заключаем, что
$
\hat{v}^\ve \in H^1(\O).
$

Дальнейшей целью будет
\begin{theorem}
Для разности  $u^\ve$ -- решения задачи (\ref{eq1}) и функции $\hat{v}^\ve$, определенной в  (\ref{eq20}), имеет место оценка
\beq\label{eq23}
\|u^\ve-\hat{v}^\ve\|_{H^1(\O)}\le C\max\{\sqrt{\ve},\sqrt{\frac{\delta}{\ve}}\}\|f\|_{L^2(\O)},
\eeq
где константа $C$ того же типа, что и в  (\ref{eq7}).
\end{theorem}

Параграфы \S 4-5 посвящены доказательству теоремы 3.1.

Для задачи Неймана с $\ve$-периодической матрицей $a^\ve(x){=}a(x/\ve)$
$H^1$-оценка порядка $O(\sqrt{\ve})$ впервые получена в \cite{PZh1}, как для скалярного уравнения так и для системы теории упругости.
В настоящей работе предложенный в \cite{PZh1} метод адаптируется для многомасштабных задач c  учетом специфики скалярного случая. Общий подход, применимый для векторных многомасштабных задач, развивался  
в \cite{PMS}.

\section{Оценка невязки}
 
Справедливо интегральное тождество
\beq\label{eq24}
\ilo (a^\ve\nabla \hat{v}^\ve - a^\ve\nabla u^\ve)\cdot\vp dx = \ilo (a^\ve\nabla \hat{v}^\ve - a^0\nabla u)\cdot\nabla\vp dx,\quad \vp\in H^1(\O),
\eeq
поскольку
$$
\ilo a^\ve\nabla u^\ve\cdot\nabla\vp dx = \ilo f\vp dx = \ilo a^0\nabla u\cdot\nabla\vp dx,\ \quad  \ \vp\in H^1(\O),
$$
в силу  уравнений (\ref{eq1}) и (\ref{eq4}).
Введем обозначение для разности 
 потоков
\beq\label{eq25}
\hat{R}_{\ve}\equiv a^\ve\nabla \hat{v}^\ve - a^0\nabla u.
\eeq
Тогда в силу (\ref{eq24}) 
\beq\label{eq26}
\ilo \hat{R}_{\ve}\cdot\nabla \vp dx= \ilo (a^\ve\nabla\hat{v}^\ve-a^0\nabla u)\cdot\nabla\vp dx,\ \vp\in H^1(\O).
\eeq
Предположим, что 
доказана оценка
\beq\label{eq260}
\ilo \hat{R}_{\ve}\cdot\nabla \vp dx
\le C\max\{\sqrt{\ve},\sqrt{\frac{\delta}{\ve}}\}\|f\|_{L^2(\O)}\|\nabla \vp\|_{L^2(\O)}.
\eeq
Тогда, полагая в ней $\vp=\hat{v}^\ve -  u^\ve$ и оценивая форму (\ref{eq26}) снизу 
в силу равенства (\ref{eq24})
и эллиптичности матрицы $a^\ve$, можно получить 
\beq\label{eq261}
\|\nabla(u^\ve-\hat{v}^\ve)\|_{L^2(\O)}\le C\max\{\sqrt{\ve},\sqrt{\frac{\delta}{\ve}}\}\|f\|_{L^2(\O)},
\eeq
что является частью оценки (\ref{eq23}).
Этим мотивирован наш интерес к оценке (\ref{eq260}).

Прежде чем оценивать сверху форму (\ref{eq26}), попробуем оценить  аналогичную форму,
в которой 
вместо $\hat{R}_{\ve}$ стоит 
разность
$$
R_\ve=a^\ve\nabla v^\ve-a^0\nabla u
$$
с    обычным первым приближением $v^\ve$, определенным в (\ref{eq18}).

Запишем  
$$
\nabla v^\ve (x)= (I+\nabla_y N(y))\nabla u (x) +\nabla_z  M(y,z)(I+\nabla_y N(y))\nabla u(x)+
$$
$$
+\ve \nabla^2 u(x)N(y)+\frac{\delta}{\ve} \nabla_y [M(y,z)\cdot(I+\nabla_y N(y))\nabla u(x)]+
$$
$$
+\delta (I+\nabla_yN(y))\nabla^2 u(x)M(y,z).
$$
Тогда
$$
R_{\ve} = a^\ve\nabla v^\ve - a^0\nabla u=
$$
\beq\label{eq27}\ds{
=a(y,z)(I+\nabla_z M(y,z)) (I+\nabla_y N(y))\nabla u(x) - \hat{a}(y)(I+\nabla_y N(y))\nabla u(x)+}
\atop\ds{
\hat{a}(y)(I+\nabla_y N(y))\nabla u(x) - a^0\nabla u(x)+r^1_{\ve},}
\eeq
где
$$
r^1_{\ve}=a(y,z)\big(\ve \nabla^2 u(x) N(y)+\f{\delta}{\ve} \nabla_y [M(y,z)(I+\nabla_y N(y))\nabla u(x)]+
$$
\beq\label{eq28}
+\ve (I+\nabla_y N(y))\nabla^2 u(x)M(y,z)\big),\ y=\frac{x}{\ve},\ z=\frac{x}{\delta}.
\eeq
Используя векторы
(\ref{eq10}),  можно записать вторую 
 разность в (\ref{eq27}) как
\[\ds{
\hat{a}(y)(I+\nabla_y N(y))\nabla u(x) - a^0\nabla u(x)=
[\hat{a}(y)(I+\nabla_y N(y))e^j-a^0e^j]\f{\pa u(x)}{\pa x_j}=}
\atop\ds{
[\hat{a}(y)(e^j+\nabla_y N_j(y))-a^0e^j]\f{\pa u(x)}{\pa x_j}\stackrel{(\ref{eq10})_2}=
g^j(y)\f{\pa u(x)}{\pa x_j}
}
\]
и первую разность в (\ref{eq27}) как
\[\ds{
a(y,z)(I+\nabla_z M(y,z)) (I+\nabla_y N(y))\nabla u(x) - \hat{a}(y)(I+\nabla_y N(y))\nabla u(x)=}
\atop\ds{
[a(y,z)(I+\nabla_z M(y,z))- \hat{a}(y)](\nabla u(x)+\nabla_y N(y)\nabla u(x) )
=}
\]
\[
\ds{
[a(y,z)(I+\nabla_z M(y,z))- \hat{a}(y)]e^j(\f{\pa u(x)}{\pa x_j}+
\f{\pa }{\pa y_j}(N(y)\cdot  \nabla u(x) )=
\atop\ds[a(y,z)(e^j+\nabla_z M_j(y,z))- \hat{a}(y)e^j]\zeta_j(x,y)
)
\stackrel{(\ref{eq10})_1}=
 p^j(y,z)\zeta_j(x,y),
}
\]
где
 \[
 \zeta_j(x,y)= \f{\pa u(x)}{\pa x_j}+\f{\pa}{\pa y_j} (N(y)\cdot\nabla u(x)).
 \]
 В итоге равенство (\ref{eq27}) представляется в  виде суммы
\beq\label{eq29}
R_{\ve} (x) = \zeta_j(x,y) p^j(y,z) +  \f{\pa u(x)}{\pa x_j}g^j(y) +r^1_\ve (x),\ y=\f{x}{\ve},\ z=\f{x}{\delta},
\eeq
где 
 первые слагаемые в силу леммы 2.2 можно записать как\\
\beq\label{eq30}
r^2_{\ve}\equiv g^j(y)\f{\pa u(x)}{\pa x_j}= \ve \Div \l(G^j(y)\f{\pa u(x)}{\pa x_j}\r) - \ve G^j(y)\nabla\l(\f{\pa u(x)}{\pa x_j}\r),
\eeq
\beq\label{eq31}
r^3_{\ve} \equiv p^j(y,z)\zeta_j(x,y) = \delta \Div ( P^j (y,z)\zeta_j(x,y)) -\f{\delta}{\ve}\Div_y (P^j (y,z)\zeta_j) - \delta P^j(y,z)\nabla_x \zeta_j,\atop
\zeta_j(x,y)= \f{\pa u(x)}{\pa x_j}+\f{\pa}{\pa y_j} (N(y)\cdot\nabla u(x)).
\eeq

Из (\ref{eq27}) -- (\ref{eq31}) получаем равенство
\beq\label{eq32}
\ilo R_{\ve}\cdot\nabla\vp dx= \sum^3_{i=1}\ilo r^i_{\ve}\cdot\nabla\vp dx.
\eeq

Введем гладкую функцию $\theta^\ve(x)$, сосредоточенную в $\ve$-окрестности границы $\Gamma=\pa\Omega$, такую что
\beq\label{eq320}
supp\,\theta^\ve\subset \Gamma_\ve,\quad
\theta^\ve (x)|_{\Gamma}=1,\quad
 0\le \theta^\ve (x)\le 1,\quad
|\ve\nabla \theta^\ve (x)|\le c.
\eeq
Используя  срезающую функцию $\theta^\ve (x)$, в (\ref{eq32}) можно заменить слагаемые $r^2_\ve$ и $r^3_\ve$ из (\ref{eq30}) и (\ref{eq31}) на их аналоги со срезающей функцией
\beq\label{eq321}\ds{
\tilde{r}^2_{\ve}=\ve\Div \left(\theta^\ve (x) G^j(y)\frac{\pa u(x)}{\pa x_j}\right)-\ve G^j(y)\nabla_x\left(\frac{\pa u(x)}{\pa x_j}\right),
}
\atop\ds{
\tilde{r}^3_{\ve}=\delta \Div(\theta^\ve (x)P^j(y,z)\zeta_j(x,y))- \f{\delta}{\ve}\Div_y (P^j (y,z)\zeta_j(x,y))- \delta P^j(y,z)\nabla_x \zeta_j(x,y).}
\eeq
Действительно, функция $(1-\theta^\ve (x))$ обращается в нуль в окрестности границы $\pa\O$
и 
вектор $\Div \l((1-\theta^\ve (x)) G^j(\f{x}{\ve})\f{\pa u(x)}{\pa x_j}\r)$ обладает свойством соленоидальности:
$$
\ilo \Div\l[(1-\theta^\ve(x))G^j(\frac{x}{\ve})\frac{\pa u(x)}{\pa x_j}\r]\cdot\nabla\vp(x) dx=
$$
\beq\label{eq33}
-\ilo \l[(1-\theta^\ve(x))G^j(\frac{x}{\ve})\cdot\nabla^2\vp(x)\r]\frac{\pa u(x)}{\pa x_j} dx=0\quad\vp\in C^\infty(\RR^d)
\eeq
(без суммирования по $j$). Равенство нулю в (\ref{eq33}) имеет место, поскольку  $G^j(\frac{x}{\ve})\cdot\nabla^2\vp(x)=0$  поточечно  в силу кососимметричности матрицы $G^j$.
Аналогичным свойством соленоидальности  обладает вектор $\Div\l[((1-\theta^\ve(x)))P^j(\frac{x}{\ve},\f{x}{\delta})\zeta_j(x,y)\r]$:
\beq\label{eq34}
\ilo\Div\l[((1-\theta^\ve(x)))P^j(\frac{x}{\ve},\f{x}{\delta})\zeta_j(x,y)\r]\cdot\nabla\vp dx =0\quad \forall\vp\in C^\infty(\O),
\eeq
где проявляется 
кососимметричность  матрицы  $P^j(y,z)$.
 
Таким образом, учитывая (\ref{eq34}) и (\ref{eq33}), из (\ref{eq32}) получаем
\beq\label{eq35}\ds{
\ilo {R}_{\ve}\cdot\nabla \vp dx = \ilo (r^1_{\ve} +\tilde{r}^2_{\ve}+\tilde{r}^3_{\ve})\cdot \nabla\vp dx\le
}
\atop\ds{
\le (\|r^1_{\ve}\|_{L^2(\O)} +\|\tilde{r}^2_{\ve}\|_{L^2(\O)}+\|\tilde{r}^3_{\ve}\|_{L^2(\O)})\|\nabla \vp \|_{L^2(\O)},\ \forall\vp\in C^\infty(\bar{\O}).}
\eeq
  Структура функций $r^1_{\ve}$, $\tilde{r}^2_\ve$ и $\tilde{r}^3_\ve$ такова, что из (\ref{eq35}) непосредственно не следует интересующая нас оценка (\ref{eq260}). Далее покажем, как обойти это препятствие.

\section{Доказательство $H^1$-оценки}

\textbf{1.}
Рассмотрим наряду с обычным первым приближением (\ref{eq18}) смещенное первое приближение:
$$
v^\ve_{\sigma} (x)= u(x')+\ve N(y')\cdot\nabla u (x')+\delta M(y',z)\cdot (I+\nabla_y N(y'))\nabla u(x'),\atop
x'=x-\delta\sigma,\ y'=\f{x}{\ve}-\f{\delta}{\ve}\sigma,\ z=\f{x}{\delta}.
$$
Ему будет соответствовать выражение
\beq\label{eq36}
R_{\ve,\sigma} (x)=a(\f{x}{\ve}-\f{\delta}{\ve}\sigma,\f{x}{\delta})\nabla v^\ve_\sigma (x)-a^0\nabla u (x-\delta\sigma)
\eeq
и аналогичная (\ref{eq35})   оценка
$$
\ilo R_{\ve,\sigma}\cdot\nabla \vp dx =  \ilo (r^1_{\ve,\sigma}+\tilde{r}^2_{\ve,\sigma}+\tilde{r}^2_{\ve,\sigma})\cdot\nabla\vp dx\le
$$
\beq\label{eq37}
\le (\|r^1_{\ve,\sigma}\|_{L^2(\O)}+\|\tilde{r}^2_{\ve,\sigma}\|_{L^2(\O)}+\|\tilde{r}^3_{\ve,\sigma}\|_{L^2(\O)})
\|\nabla\vp\|_{L^2(\O)}.
\eeq
Проинтегрировав соотношение (\ref{eq37}) по  $\sigma\in Z$, получаем
$$
\int\limits_Z\ilo R_{\ve,\sigma}\cdot\nabla \vp dx d\sigma \le
$$
\beq\label{eq38}
\le \|\nabla\vp\|_{L^2(\O)}\sum\limits_{j}\gamma_j\left(\ilo\int\limits_Z\left|b_j \left(\frac{x}{\ve}-\frac{\delta}{\ve}\sigma,\frac{x}{\delta}\right)U_j (x-\delta\sigma)\right|^2  d\sigma dx\right)^\frac12.
\eeq
Здесь коэффициент $\gamma_{j}$ равен  $\ve$, $\frac{\delta}{\ve}$, $\delta$ или $\ve^0$. Элементы $b_{j}(y,z)$ сформированы из  $N_j(y)$,  $M_j(y,z)$,  $P^j(y,z)$,  $G^j(y)$ и их  производных, а  $U_{j}(x)$ состоят из градиентов $\nabla u$, $\nabla^2 u$,  умноженных, быть может, на ограниченные функции $\theta^\ve (x)$ или $\ve \nabla\theta^\ve (x)$. 
Оценим сверху интегралы из суммы (\ref{eq38}).

Например, выражение (см. (\ref{eq321}))
$$
\tilde{r}^3_{\ve,\sigma}=\delta\Div \l[\theta^\ve (x')P^j(y',z)\zeta_j(x',y')\r]-\f{\delta}{\ve} \Div_y\l(P^j(y',z)\zeta_j(x',y')\r)-\delta P^j(y',z)\nabla_x\zeta_j(x',y'),
$$
где
$$
\zeta_j(x',y')=\f{\pa u(x')}{\pa x_j} +\f{\pa}{\pa y_j}(N(y')\cdot\nabla u(x')),\ x'=x-\delta\sigma,\ y'=\f{x}{\ve}-\f{\delta}{\ve}\sigma,\ z=\f{x}{\delta},
$$
дает вклад в сумму (\ref{eq38}) в виде интегралов
\beq\label{eq39}
\ds{\ilo\int\limits_Z\left|p^j(y',z)\theta^\ve (x')\zeta_j(x',y')\right|^2 d\sigma dx,\ \ilo\int\limits_Z\left|P^j(y',z)(\delta\nabla\theta^\ve (x'))\zeta_j(x',y')\right|^2 d\sigma dx,}
\eeq
\beq\label{eq40}
\ds{ \ilo\int\limits_Z \left|\Div_y(P^j(y',z)\zeta_j(x',y'))\right|^2 d\sigma dx,\ \ilo\int\limits_Z \left|P^j(y',z)\nabla_x\zeta_j(x',y')\right|^2 d\sigma dx.}
\eeq

Интегралы  из (\ref{eq39}) не  содержат малого множителя  типа $\ve$, $\delta$  в сумме (\ref{eq38}), но они сами малы, так как благодаря срезающим функциям их можно оценить, используя неравенство для следа
\beq\label{eq41}
\|\vp\|^2_{L^2(\Gamma_{\ve})}\le c \ve  \|\vp\|_{L^2(\O)}\|\nabla \vp \|_{L^2(\O)},\ c=const(d,\O).
\eeq
Например,  для одного из этих  интегралов  имеем
$$
\ilo\int\limits_{Z}\left|P^j(y',z)(\delta\nabla\theta^\ve (x'))\zeta_j(x',y')\right|^2  d\sigma dx\stackrel{(\ref{eq320})}\le
$$
$$
\le \int\limits_{\Gamma_{\ve}}\int\limits_{Z}\left|P^j(y',z)\zeta_j(x',y')\right|^2  d\sigma dx\le  \int\limits_{\Gamma_{2\ve}}\int\limits_{Z}\left|P^j(y,z+\sigma)\zeta_j(x,y)\right|^2 d\sigma  dx\le
$$
$$
\le c \int\limits_{\Gamma_{2\ve}}(|u(x)|^2 +|\nabla u(x)|^2) dx
\stackrel{(\ref{eq41})}\le  c_2\ve(\|u\|_{L^2(\O)}\|\nabla u\|_{L^2(\O)}+ \|\nabla u\|_{L^2(\O)}\|\nabla^2 u\|_{L^2(\O)})\le C\ve \|f\|^2_{L^2(\O)}.
$$
Здесь на втором и третьем шаге использованы те же соображения, что  при выводе леммы 3.1 -- замена переменной и интегрирование по дополнительному параметру $\sigma$.

Оценим один из интегралов  типа (\ref{eq40}):
$$
\ilo\int\limits_Z \left|\Div_y(P^j(y',z)\zeta_j(x',y'))\right|^2 d\sigma dx \le
$$
  $$ \int\limits_{\O_{\ve}}\int\limits_Z \l|(\Div_y P^j(y,z+\sigma))\zeta_j(x,y)+P^j(y,z+\sigma)\nabla_y \zeta_j(x,y)\r|^2 d\sigma dx\le
$$
$$
c\int\limits_{\O_{\ve}}\int\limits_Z \l|\Div_y P^j(y,z+\sigma)|^2|\nabla u(x)\r|^2  d\sigma dx + \int\limits_{\O_\ve}\int\limits_Z \l|P^j(y,z+\sigma)\nabla^2 N (y)\nabla u(x)\r|^2 d\sigma dx\le 
$$
\[
c_1\sup_y\|\nabla_y P^j(y,\cdot)\|^2_{L^2(Z)}\|\nabla u\|^2_{L^2(\O)}+
\sup_y\| P^j(y,\cdot)\|^2_{L^2(Z)} \int\limits_{\Omega_1 } |\nabla^2 N (\f{x}{\ve})\nabla u(x)|^2 dx\le C\|f\|^2_{L^2(\O)}.
\]
где  использованы свойства  функций $N_j$, $P^j$ и $u$ (см. лемму 2.1 
и оценку (\ref{eq5})),
оценки типа (\ref{eq210}) и (\ref{eq230}), а также свойства оператора продолжения (\ref{eq6}).

Все оставшиеся интегралы в (\ref{eq39}) и (\ref{eq40}) можно оценить, применяя аналогичные соображения. То же верно для других интегралов в сумме (\ref{eq38}), возникающих из представления типа (\ref{eq35}).
В итоге 
 имеем оценку
\beq\label{eq42}
\left|\int\limits_Z\ilo R_{\ve,\sigma}\cdot\nabla\vp dx d\sigma\right|\le C\max\{\sqrt{\ve},\sqrt{\frac{\delta}{\ve}}\}\|f\|_{L^2(\O)}\|\nabla\vp\|_{L^2(\O)}.
\eeq

\textbf{2.}
Теперь попробуем заменить в (\ref{eq42}) $R_{\ve,\sigma}$ на выражение $\hat{R}_{\ve}$, определенное в (\ref{eq25}).
Справедливо следующее неравенство
\beq\label{eq43}
\|\int\limits_Z R_{\ve,\sigma}  d\sigma - \hat {R}_{\ve}\|_{L^2(\O)}\le C\max\{\ve,\frac{\delta}{\ve}\} \|f\|_{L^2(\O)}.
\eeq
Действительно, запишем  $\int\limits_Z R_{\ve,\sigma}  d\sigma$ подробно
\[
\int\limits_Z R_{\ve,\sigma} d\sigma=\int\limits_Z a(y-\f{\delta}{\ve}\sigma,z)\bigg(\nabla u(x-\delta\sigma)+\ve\nabla [N(y-\f{\delta}{\ve}\sigma)\cdot\nabla u(x-\delta\sigma)]+
\] 
$$
+\delta\nabla[M(y-\f{\delta}{\ve}\sigma,z)\cdot(I+\nabla N(y-\f{\delta}{\ve}\sigma))\nabla u(x-\delta\sigma)]\bigg)d\sigma -a^0 \int\limits_Z\nabla u(x-\delta\sigma)d\sigma,\ y=\frac{x}{\ve},\ z=\frac{x}{\delta}.
$$
 Заменим матрицу $a\l(\frac{x}{\ve}-\frac{\delta}{\ve}\sigma,\frac{\delta}{\ve}\r)$  на $a\left(\frac{x}{\ve},\frac{x}{\delta}\right)$ с погрешностью $ O(\frac{\delta}{\ve})$ в силу свойства липшицевости, в результате чего, используя обозначение для корректора из  
 (\ref{eq1.3.3.44}), имеем
$$
\int\limits_Z R_{\ve,\sigma}(x) d\sigma=
$$
\beq\label{eq45}
=a^\ve(x)[\int\limits_Z \nabla u(x-\delta\sigma)d\sigma+\nabla K_{\ve}(x)]-
a^0  \int\limits_Z \nabla u(x-\delta\sigma) d\sigma+ O(\frac{\delta}{\ve}).
\eeq
Заметим, что
$
\int\limits_Z \varphi(x-\delta\sigma)d\sigma{ = }(\varphi)_{\delta} (x)
$
есть  сглаживание по Стеклову функции $\varphi(x)$. Известно, что $(\varphi)_{\delta} (x)$
аппроксимирует саму функцию $\varphi(x)$ с оценкой
$$
\|(\varphi)_{\delta}-\varphi\|_{L^2(\O)}\le c\delta\|\nabla \varphi\|_{L^2(\O)}, c=const(d).
$$
Поэтому в (\ref{eq45}) с учетом оценки (\ref{eq5})   опускаем
сглаживание по Стеклову для $\nabla u$, после чего
приходим к соотношению 
$$
\int\limits_Z R_{\ve,\sigma}(x) d\sigma=
a^\ve(x)(\nabla u(x)+\nabla K_{\ve}(x))-a^0\nabla u(x)+O(\max\{\delta,\frac{\delta}{\ve}\})\stackrel{(\ref{eq1.3.3.44})}=$$
$$
a^\ve(x)\hat{v}^{\ve}(x)-a^0\nabla u(x)
+O(\max\{\ve,\frac{\delta}{\ve}\})\stackrel{(\ref{eq25})}=
\hat{R}_{\ve}
+O(\max\{\ve,\frac{\delta}{\ve}\}),
$$
то есть оценка (\ref{eq43}) получена.

В итоге из (\ref{eq42}), (\ref{eq43}) имеем
\beq\label{eq46}
\ilo \hat{R}_{\ve}\cdot\nabla \vp dx\le C \max\{\sqrt{\ve},\sqrt{\frac{\delta}{\ve}}\}\|f\|_{L^2(\O)}\|\nabla\vp\|_{L^2(\O).}
\eeq
Таким образом,  установлены соотношение (\ref{eq260}) и как следствие оценка (\ref{eq261}).

\textbf{3.} Теперь по неравенству Пуанкаре  записываем оценку
\beq\label{eq47}
\|\hat{v}^\ve-u^\ve\|^2_{L^2(\O)}\le C_P\l(\|\nabla(\hat{v}^\ve-u^\ve)\|^2_{L^2(\O)}
+\l(\ilo (\hat{v}^\ve-u^\ve) dx\r)^2\r).
\eeq
Учитывая, что $\ilo u^\ve dx=\ilo u dx=0$, имеем
$$
\ilo (\hat{v}^\ve-u^\ve) dx = \ilo K_{\ve} dx.
$$
В \S 3 фактически было доказано свойство 
\[
\| K_{\ve}\|_{L^2(\O)}\le C\max\{\ve,\delta\}\|f\|_{L^2(\O)}.
\]
Таким образом, из (\ref{eq47}) и (\ref{eq261}) вытекает оценка
\[
\|\hat{v}^\ve-u^\ve\|^2_{L^2(\O)}\le C\max\{\ve,\frac{\delta}{\ve}\}\|f\|^2_{L^2(\O)}, 
\]
и в итоге  получаем $H^1$-оценку (\ref{eq23}). Теорема 4.1 доказана.

\textbf{4.}
Оценку (\ref{eq46}) можно записать в несколько уточненном виде. Для этого  в сумме (\ref{eq38}) слагаемые с множителями $\ve\nabla\theta^\ve (x)$ и $\theta^\ve (x)$ не надо огрублять так, как раньше. В результате получается оценка
\beq\label{eq48}
\left|\ilo \hat{R}_{\ve}\cdot\nabla\vp dx\right|\le C\|f\|_{L^2(\O)}(\sqrt{\ve}\|\nabla\vp\|_{L^2(\Gamma_{2\ve})}+\max\{{\ve},{\frac{\delta}{\ve}}\}\|\nabla\vp\|_{L^2(\O)}).
\eeq

При выводе (\ref{eq46})  важную роль играли следующие свойства нулевого приближения $u$:\\
$\quad{}\quad{}\quad$ 1) $u\in H^2(\O)$ с оценкой $\|u\|_{H^2(\O)}\le c_0\|f\|_{L^2(\O)}$;\\
$\quad{}\quad{}\quad$  2) неравенство для градиента $\nabla u$ в окрестности границы $\pa\O$
\beq\label{eq49}
\int\limits_{\Gamma_{\ve}}|\nabla u|^2 dx\le c\ve \|f\|^2_{L^2(\O)}.
\eeq
Второе свойство есть следствие первого.
Отметим, что  свойство (\ref{eq49}) справедливо и для градиента исходной задачи $\nabla u^\ve$, хотя $\nabla u^\ve$ необязательно принадлежит $H^2(\O)$.

\begin{lemma}
Пусть $u^\ve$ -- решение задачи (\ref{eq1}). Тогда выполнена оценка
$$
\int\limits_{\Gamma_{\ve}}|\nabla u^\ve|^2  dx \le C\max\{\ve,\frac{\delta}{\ve}\} \|f\|^2_{L^2(\O)},
$$
где $C$ того же типа, что и в (\ref{eq7}).
\end{lemma}
\doc
Запишем равенство 
$$ u^\ve = (u^\ve-\hat{v}^\ve)+\hat{v}^\ve = (u^\ve-\hat{v}^\ve)+u +K_{\ve},$$
где $\hat{v}^\ve$ есть $H^1$-приближение к решению $u^\ve$ (см. (\ref{eq1.3.3.44})).
Отсюда
$$
\nabla u^\ve =\nabla (u^\ve -\hat{v}^\ve) +\nabla u +\nabla K_{\ve},
$$
$$
\int\limits_{\Gamma_{\ve}}|\nabla u^\ve|^2 dx \le \int\limits_{\Gamma_{\ve}} |\nabla (u^\ve-\hat{v}^\ve)|^2 dx + \int\limits_{\Gamma_{\ve}} |\nabla u|^2 dx +\int\limits_{\Gamma_{\ve}} |\nabla  K_{\ve}|^2 dx= I_1+I_2+I_3.
$$
Оценим каждый из интегралов $I_1$, $I_2$ и $I_3$:
$$
I_1\le \|u^\ve - \hat{v}^\ve\|^2_{H^1(\O)}\stackrel{(\ref{eq23})}\le C\max\{\ve,\frac{\delta}{\ve}\} \|f\|^2_{L^2(\O)};
$$
$$
I_2 \stackrel{(\ref{eq49})}\le c\ve \|f\|^2_{L^2(\O)};
$$
$$
I_3\le \int\limits_{\Gamma_{2\ve}} (
|\nabla u|^2+\max\{\ve^2,\l(\frac{\delta}{\ve}\r)^2\} |\nabla^2 u|^2) dx \le
$$
$$
\le c(\ve\|f\|^2_{L^2(\O)}+\max\{\ve^2,\l(\frac{\delta}{\ve}\r)^2,\delta^2\} \|f\|^2_{L^2(\O)},
$$
что  устанавливается  с помощью рассуждений  из \S3 и неравенства для следа.

\section{Вывод $L^2$-оценки}
Для произвольной  функции $\Phi\in \tilde{L}^2(\O)$ рассмотрим  краевую задачу Неймана
\beq\label{eq50}
\vp^\ve\in \tilde{H}^1(\Omega),\quad -\Div [a^\ve_T (x)\nabla \vp^\ve(x)]=\Phi,
\eeq
где $a^\ve_T (x)$ -- транспонированная к (\ref{eq2})  матрица.

Для решения $\vp^\ve$ выполнены энергетическая оценка 
\beq\label{eq51}
\|\vp^\ve\|_{H^1(\O)}\le c_0\|\Phi\|_{L^2(\O)}
\eeq
и оценка для градиента вблизи границы (см. лемму 5.1)
\beq\label{eq52}
\|\nabla\vp^\ve\|_{L^2(\Gamma_{2\ve})}\le C\max\{\sqrt{\ve},\sqrt{\frac{\delta}{\ve}}\}\|\Phi\|_{L^2(\O)}.
\eeq

Исходя из интегрального тождества для задачи  (\ref{eq50}) с пробной функцией 
$$w^\ve = \hat{v}^\ve-u^\ve$$ имеем
$$
\ilo \Phi w^\ve dx=
\ilo a^\ve_T\nabla\vp^\ve\cdot\nabla w^\ve dx =  \ilo \nabla \vp^\ve\cdot a^\ve\nabla w^\ve dx
\stackrel{(\ref{eq24})}= \ilo \hat{R}_{\ve}\cdot\nabla \vp^\ve dx,
$$
где в силу (\ref{eq48}), (\ref{eq51}), (\ref{eq52})
$$
\ilo \hat{R}_{\ve}\cdot\nabla\vp^\ve dx\le C\|f\|_{L^2(\O)}(\sqrt{\ve}\|\nabla\vp^\ve\|_{L^2(\Gamma_{2\ve})}+\max\{\ve,\frac{\delta}{\ve}\}\|\nabla\vp^\ve \|_{L^2(\O)})\le
$$
$$
\le C\max\{\ve,{\frac{\delta}{\ve}}\}\|\Phi\|_{L^2(\O)}\|f\|_{L^2(\O)}.
$$
Отсюда, взяв $\Phi=w^\ve-\int_\Omega w^\ve\,dx$, получаем
\[
\|\Phi\|_{L^2(\O)}\le C\max\{\ve,{\frac{\delta}{\ve}}\}\|f\|_{L^2(\O)}.
\]
В итоге, поскольку
$
\int_\Omega w^\ve\,dx=\int_\Omega K_\ve\,dx=O(\ve)$, что уже показано раньше,
выводим
$$
\|\hat{v}^\ve-u^\ve\|_{L^2(\O)}\le C\max\{\ve,{\frac{\delta}{\ve}}\}\|f\|_{L^2(\O)}.
$$
Отсюда уже получается $L^2$-оценка  (\ref{eq7}), так как
$$
\|u^\ve-u\|_{L^2(\O)}=\|u^\ve-\hat{v}^\ve +K_\ve\|_{L^2(\O)}\le \|u^\ve-\hat{v}^\ve\|_{L^2(\O)} + O(\ve).
$$
Теорема 1.1 доказана.

Для задачи Неймана с $\ve$-периодической матрицей $a^\ve(x){=}a(x/\ve)$
$L^2$-оценка порядка $O({\ve})$ впервые получена в \cite{Su2}.

\section{Доказательство вспомогательных утверждений}
Приведем доказательство леммы 2.1.

Из интегрального тождества для  первой задачи на ячейке 
$$
\Div_z [a(y,z) (e^j+\nabla_z M_j(y,z))]=0
$$
легко получить равенство форм
$$
 \int\limits_Z a(y,z)(e^j+\nabla_z M_j(y,z))\cdot(e^j+\nabla_z M_j(y,z)) dz=\int\limits_{Z}a(y,z)(e^j+\nabla_z M_j (y,z))\cdot e^j dz.
$$
Оценивая левую форму снизу, а правую -- сверху
по условию (\ref{eq3})$_2$ на 
матрицу $a(y,z)$,  имеем
$$
\mu\int\limits_Z |e^j+\nabla_z M_j(y,\cdot)|^2\,dz
\le\mu^{-1}\int\limits_Z |e^j+\nabla_z M_j(y,\cdot)|\,dz
\le \mu^{-1}(\int\limits_Z |e^j+\nabla_z M_j(y,\cdot)|^2\,dz)^{\frac12},
$$
что дает 
\beq\label{eq58}
(\int\limits_Z |e^j+\nabla_z M_j(y,\cdot)|^2\,dz)^{\frac12}\le \mu^{-2}.
\eeq
Отсюда следует, что
$$
\|M_j(y,\cdot)\|_{H^1_{per}(Z)}\le c,\ c=const(\mu),
$$
и  свойство $i)$ доказано.

Докажем липшицевость функции $M_j(y,z)$ по переменной $y$. Положим:
$$
M^{h}_j(y,z)= M_j(y+h,z),\ a^h (y,z)=a(y+h,z).
$$
Для решений уравнений
$$
\Div_z [a^h (y,z)(e^j+\nabla_z M_j^{h}(y,z))]=0,\quad \Div_z [a(y,z)(e^j+\nabla_z M_j(y,z))]=0.
$$
запишем интегральные тождества  на пробной функции 
$ \varphi (z) = M_j^{h}(\cdot,z)-M_j(\cdot,z)$.
Вычитанием получаем
$$
\int\limits_{Z}(a^h (y,z)(e^j+\nabla_z M_j^{,h}(y,z))-a(y,z)(e^j+\nabla_z M_j(y,z)))\cdot\nabla\varphi(z) dz=0,
$$
$$
\int\limits_{Z}a^h\nabla_z (M_j^{h}-M_j)\cdot\nabla_z (M_j^{h}-M_j)dz=-\int\limits_Z (a^h-a)(e^j+\nabla_z M_j)\cdot \nabla_z(M_j^{h}-M_j)dz.
$$
Используя условие  (\ref{eq3})$_1$ и доказанное  уже свойство 
 ограниченности  $M(y,z)$, выводим
$$
\mu\int\limits_{Z}|\nabla_z(M_j^{h}-M_j)|^2dz\le
 c_L |h|(\int\limits_{Z}|e^j+\nabla_z M_j|^2dz)^{\frac12}(\int\limits_Z |\nabla_z (M_j^{h}-M_j)|^2 dz)^\frac12\stackrel {(\ref{eq58})}\le
$$
$$
 \mu^{-2}|h|c_L(\int\limits_Z |\nabla_z (M_j^{h}-M_j)|^2dz)^\frac12,
$$
Отсюда
$$
(\int\limits_{Z}|\nabla_z(M_j^{h}-M_j)|^2dz)^\frac12\le \mu^{-3} c|h|,
$$
и липшицевость $M_j(\cdot,z)$ как функции со значениями в $H^1_{per}(Z)$ доказана.

Докажем липшицевость матрицы $\hat{a}(y)$, где
$$
\hat{a}(y)=\langle a(y,\cdot)(I+\nabla_z M(y,\cdot))\rangle_Z,
$$
$$
\hat{a}^h(y)=\hat{a}(y+h)=\langle a(y+h,\cdot)(I+\nabla_z M(y+h,\cdot))\rangle_Z.
$$
Для разности $\hat{a}^h (y) - \hat{a}(y)$ имеем
$$
(\hat{a}^h (y) - \hat{a}(y))=\int\limits_Z (a^h(y,z)(I+\nabla_z M^{h})-a(y,z)(I+\nabla_z M))dz=
$$
$$
=\int\limits_Z [a^h(I+\nabla_z M^{h})-a(I+\nabla_z M^{h})+a(I+\nabla_z M^{h})-a(I+\nabla_z M)]dz=
$$
$$
=\int\limits_Z [(a^h-a)(I+\nabla_z M^{h})-a\nabla_z (M-M^{h})]dz.
$$
Откуда, в силу доказанных уже свойств $i)$ и $ii)$, следует
$$
|\hat{a}^h (y)-\hat{a}(y)|\le c_L|h|\int\limits_Z |I+\nabla_z M^{h}| dz +\mu^{-1}\int\limits_{Z} |\nabla_z (M^{h}-M)|dz\le c|h|,
$$
 что и требуется.
 
  Уравнение (\ref{eq9}) можно записать как
  \[
  \Div_y (\hat{a}(y)\nabla N_j(y))=-\Div_y (\hat{a}(y)e^j).
  \]
  Поскольку  $\hat{a}(y)$ липшицева и, значит, $\Div_y (\hat{a}(y)e^j)\in L^p(Y)$, к нему применима эллиптическая теория  \cite[гл. III, \S 15]{LaUr}, что обеспечивает свойство  $v)$.

Перейдем к утверждениям $iv)$ и $vi)$. Ограниченность решений обеих задач на ячейке следует из обобщенного принципа максимума \cite[гл. II, Приложение B]{KS}. Покажем,  что  ограничен также градиент  $\nabla_y N_j(y)$. Продифференцируем уравнение (\ref{eq9}), что возможно, благодаря липшицевости матрицы $\hat{a}(y)$. Имеем
$$
\Div_y\l[\hat{a}(y)\nabla_y \l(\f{\pa N_j (y)}{\pa y_k}\r)+\f{\pa\hat{a}(y)}{\pa y_k} (e^j +\nabla_y N_j(y))\r]=0, \quad k=1,\ldots d.
$$
Вводим обозначение
$v_k=\f{\pa N_j}{\pa y_k}$ и получаем систему относительно вектора $v=\{v_k\}$,
$$
\Div_y[\hat{a}(y)\nabla v_k+\f{\pa\hat{a}(y)}{\pa y_k}(e^j+ v)]=0,\quad k=1,\ldots d.
$$
Известно \cite[гл. VII, \S 2]{LaUr}, что решение $v$ такой диагональной системы принадлежит $L^\infty(Y)$. Отсюда выводим, что $\nabla_y N_j\in L^\infty(Y)$.

Лемма доказана.

{}

\end{document}